\documentclass[10pt]{amsart}
\usepackage{url,graphicx}
\usepackage[dvips]{epsfig}
\usepackage{latexsym,color}
\usepackage{amscd,amssymb,verbatim}
\newcommand{\be}{\begin{enumerate}}
\newcommand{\ee}{\end{enumerate}}

\newcommand{\bd}{{\text{\rm Bd}\,}}

\newcommand{\cn}{{\mathcal N}}
\newcommand{\cp}{{\mathcal P}}

\newcommand{\thom}{\text{\tt Hom}\,}
\newcommand{\ra}{\rightarrow}
\newcommand{\nin}{\noindent}
\newcommand{\pr}{\noindent{\bf Proof. }}

\newtheorem{thm}{Theorem}[section]
\newtheorem{df}[thm]{Definition}

\newtheorem{rem}[thm]{Remark}

\numberwithin{equation}{section}
\numberwithin{figure}{section}

\begin{document}

\title[A simple proof for folds in complexes 
of graph homomorphisms]
{A simple proof for folds on both sides in complexes 
of graph homomorphisms}

\author{Dmitry N. Kozlov}
\address{Department of Computer Science, Eidgen\"ossische Technische
Hochschule, Z\"urich, Switzerland}
\email{dkozlov@inf.ethz.ch}
\thanks {Research supported by Swiss National Science Foundation Grant PP002-102738/1}

\date\today

\begin{abstract}
 
  In this paper we study implications of folds in both parameters 
  of Lov\'asz' $\thom(-,-)$ complexes. There is an important 
  connection between the topological properties of these  
  complexes and lower bounds for chromatic numbers.
  We give a~very short and conceptual proof of the fact that 
  if $G-v$ is a~fold of $G$, then $\bd\thom(G,H)$ collapses 
  onto $\bd\thom(G-v,H)$, whereas $\thom(H,G)$ collapses 
  onto $\thom(H,G-v)$.
  
  We also give an easy inductive proof of the only nonelementary
  fact which we use for our arguments: if $\varphi$ is a~closure
  operator on $P$, then $\Delta(P)$ collapses onto $\Delta(\varphi(P))$.
  
\end{abstract}

\maketitle

\section{Introduction.}

The $\thom$ complexes were defined by Lov\'asz. It has been shown in
\cite{BK03a,BK03b,BK03c} that the algebro-topological invariants
of $\thom$´s can be used to provide lower bounds for chromatic numbers
of graphs, a~notoriously difficult problem. Such a complex
$\thom(G,H)$ depends on two parameters, which are both (not
necessarily simple) graphs.  $\thom(G,H)$ is usually not simplicial;
however, it is a regular CW complex, and its cells are products of
simplices. The topological and combinatorial properties of the $\thom$
complexes have been studied extensively in a~recent series of papers
\cite{BK03a,BK03b,BK03c,Cs,CK1,CK2,Do,Ko1}.

The spectral sequence computations in \cite{BK03c} required certain
manipulations with the first parameter graph; these specific 
manipulations are usually called {\it folds}. It was proved 
in \cite[Proposition 5.1]{BK03b} that folds in the first 
parameter yield homotopy equivalence. It was noticed in 
\cite[Lemma 3.1]{CK1},\cite{Do} that we can fold in the second
parameter if the deleted vertex is an identical twin. More 
recently, it was observed by Csorba and Dochtermann, \cite{Cs,Do}, 
that a~proof of \cite[Proposition 5.1]{BK03b} can be modified 
to cover the folds in the second parameter completely, with 
further complications arising in certain situations. 
All of the previously known proofs follow the same path 
of using the \cite[Proposition 3.2]{BK03b}, in turn based 
on the Quillen's theorem~A, \cite[p.~85]{Qu}.

In this paper we simplify and generalize the situation. We prove
that folds in the first parameter induce a~collapsing (a~sequence
of elementary collapses) on the barycentric subdivision of the 
involved complexes, whereas folds in the second parameter induce
a~collapsing on the complexes themselves.

Our proof is simple and conceptual. We also derive in an elementary
inductive way the only fact of the topological combinatorics which we
use: the existence of the collapsing induced by a~closure operator
on a~poset.

\section{Closure operators.}

For a poset $P$ we let $\Delta(P)$ denote its {\it nerve}: the
simplicial complex whose simplices are all chains of $P$. For
a~regular CW complex $X$ we let $\cp(X)$ denote its face poset, in
particular, $\Delta(\cp(X))=\bd X$. By analogy, we denote
$\cp(\Delta(P))=\bd P$.  A~cellular map $\varphi:X\ra Y$ between
regular CW complexes induces an order preserving map
$\cp(\varphi):\cp(X)\ra\cp(Y)$.  For a~simplicial complex $X$ and its
subcomplex $Y$ we say that $X$ {\it collapses onto} $Y$ if there
exists a~sequence of elementary collapses leading from $X$ to $Y$.  If
$X$ collapses onto $Y$, then $Y$ is a~strong deformation retract
of~$X$.

Recall that an order preserving map $\varphi$ from a poset $P$ to
itself is called a~{\it descending closure operator} if
$\varphi^2=\varphi$ and $\varphi(x)\leq x$, for any $x\in P$;
analogously, $\varphi$ is called an~{\it ascending closure operator}
if $\varphi^2=\varphi$ and $\varphi(x)\geq x$, for any $x\in P$. That
ascending and descending closure operators induce strong deformation
retraction is well known in topological combinatorics, see e.g.,
\cite[Corollary 10.12]{Bj}, where it is proved by using the Quillen's
theorem~A, \cite[p.\ 85] {Qu}.

As a sample, it yields an extremely short proof of the fact that the
complex of disconnected graphs is a strong deformation retract of the
order complex of the partition lattice: the ascending closure operator
takes each disconnected graph $G$ to the partition, whose blocks are
the connected components of $G$. We remark that the first complex
appeared in the work of Vassiliev on knot theory, \cite{Va}, whereas
the second complex encodes the geometry of the braid arrangement by
means of the Goresky-MacPherson theorem, see~\cite{GM}.

In this paper we give a~short and self-contained inductive proof of
the fact that $\Delta(P)$ collapses onto $\Delta(\varphi(P))$.

\begin{thm} \label{thm1}
Let $P$ be a poset, and let $\varphi$ be a descending closure operator,
then $\Delta(P)$ collapses onto $\Delta(\varphi(P))$.
By symmetry the same is true for an ascending closure operator.
\end{thm}

\pr We use induction on $|P|-|\varphi(P)|$. If $|P|=|\varphi(P)|$,
then $\varphi$ is the identity map and the statement is obvious.
Assume that $P\setminus\varphi(P)\neq\emptyset$ and let $x\in P$ be 
one of the minimal elements of $P\setminus\varphi(P)$. 

\begin{figure}[hbt]
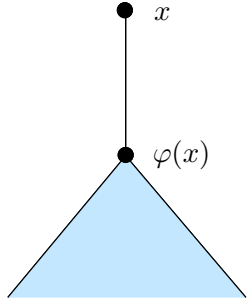

\begin{center}
  \begin{picture}(0,0)%
    \includegraphics{phi.pstex}%
  \end{picture}%
  \input{phi.pstex_t}%
  
\end{center}
\caption{$P_{<x}=P_{\leq\varphi(x)}$.}
\label{fig:phi}
\end{figure}

Since $\varphi$ fixes each element in $P_{<x}$, $\varphi(x)<x$, and
$\varphi$ is order preserving, we see that $P_{<x}$ has $\varphi(x)$
as a~maximal element, see Figure~\ref{fig:phi}. Thus the link of $x$
in $\Delta(P)$ is
$\Delta(P_{>x})*\Delta(P_{<x})=\Delta(P_{>x})*\Delta(P_{<\varphi(x)})*{\varphi(x)}$,
in particular, it is a cone with apex~$\varphi(x)$.

Let $\sigma_1,\dots,\sigma_t$ be the simplices of
$\Delta(P_{>x})*\Delta(P_{<\varphi(x)})$ ordered so that 
the dimension is weakly decreasing. Then
$$(\sigma_1\cup\{x\},\sigma_1\cup\{x,\varphi(x)\}),\dots,
(\sigma_t\cup\{x\},\sigma_t\cup\{x,\varphi(x)\})$$
is a sequence of
elementary collapses leading from $\Delta(P)$ to
$\Delta(P\setminus\{x\})$. Since $\varphi$ restricted to
$P\setminus\{x\}$ is again a descending closure operator,
$\Delta(P\setminus\{x\})$ collapses onto
$\Delta(\varphi(P\setminus\{x\}))=\Delta(\varphi(P))$ by the induction
assumption.  \qed

\begin{rem}
  There is a direct way to describe the elementary collapses in the
  proof of the Theorem \ref{thm1}. For $x\in\bd
  P\setminus\bd\varphi(P)$, $x=(x_1<\dots<x_k)$, let $1\leq i\leq k$
  be the minimal possible index, such that $x_i\notin\varphi(P)$. Then 
  either $i=1$ or $x_{i-1}\in\varphi(P)$. If $\varphi(x_i)=x_{i-1}$,
  then match $(x\setminus\{x_{i-1}\},x)$, otherwise match
  $(x,x\cup\varphi(x_i))$; the latter is possible since either $i=1$,
  or $x_i>x_{i-1}$ and $\varphi(x_i)\neq x_{i-1}$ imply
  $\varphi(x_i)>\varphi(x_{i-1})=x_{i-1}$. It is now easy to see that
  this is an acyclic matching and thus alternatively derive the result
  by using discrete Morse theory, see~\cite{Fo}.
\end{rem}

\section{$\thom$ complexes and folds.}

Let $G$ be a graph, whose set of vertices is denoted by $V(G)$, and
set of edges - by $E(G)$. Let $v$ be a~vertex of $G$, and let $G-v$
denote the graph which is obtained from $G$ by deleting the vertex
$v$ and all adjacent to $v$ edges. Let $\cn(v)$ denote the set of
neighbors of $v$, i.e., $\cn(v)=\{x\in V(G)\,|\,(v,x)\in E(G)\}$.

\begin{df} $G-v$ is called a {\bf fold} of $G$ if there exists
$u\in V(G)$, $u\neq v$, such that $\cn(u)\supseteq\cn(v)$.
\end{df}

Let $G-v$ be a fold of $G$. We let $i:G-v\hookrightarrow G$ denote the
inclusion homomorphism, and let $f:G\ra G-v$ denote the folding
homomorphism defined by $f(v)=u$ and $f(x)=x$, for $x\neq v$.

The following construction, which is due to Lov\'asz has been in the
center of latest advances in the area of topological obstructions to
graph colorings, see \cite{BK03a,BK03b,BK03c,Cs,CK1,CK2,Do,Ko1}.

Let $\Delta^{V(H)}$ be a~simplex whose set of vertices is $V(H)$. Let
$C(G,H)$ denote the direct product $\prod_{x\in V(G)}
\Delta^{V(H)}$, i.e., the copies of $\Delta^{V(H)}$ are indexed by
vertices of~$G$.

\begin{df} \label{dfhom} \cite[Definition 2.1.5]{Ko1}. 
$\thom(G,H)$ is the subcomplex of $C(T,G)$ defined by the following
condition: $\sigma=\prod_{x\in V(G)}\sigma_x\in\thom(G,H)$ if and only
if for any $x,y\in V(G)$, if $(x,y)\in E(G)$, then
$(\sigma_x,\sigma_y)$ is a~complete bipartite subgraph of~$H$, i.e., 
$\sigma_x\times\sigma_y\subseteq E(H)$.
\end{df}

Note, that $\thom(G,H)$ is a polyhedral complex whose cells are
indexed by all functions $\eta:V(G)\rightarrow
2^{V(H)}\setminus\{\emptyset\}$, such that if $(x,y)\in E(G)$, then
$\eta(x)\times\eta(y)\subseteq E(H)$.  The closure of a~cell $\eta$
consists of all cells $\tilde\eta$, satisfying
$\tilde\eta(v)\subseteq\eta(v)$, for all $v\in V(G)$.

It was noticed in \cite{BK03b} that $\thom(G,-)$ is a~covariant
functor, while $\thom(-,G)$ is a~contravariant functor from the
category of graphs and graph homomorphisms, {\bf Graphs} to {\bf Top}.
For a~graph homomorphism $\varphi:G\ra G'$ the topological maps
induced by composition are denoted as
$\varphi^H:\thom(H,G)\ra\thom(H,G')$ and
$\varphi_H:\thom(G',H)\ra\thom(G,H)$, this is the notation introduced
in \cite{BK03b}. Observe that $\varphi^H$ is cellular on the mentioned
complexes, whereas $\varphi_H$ is cellular on their first barycentric
subdivisions.

\begin{thm}\label{thm2}
  Let $G-v$ be a fold of $G$ and let $H$ be some graph. Then
  $\bd\thom(G,H)$ collapses onto $\bd\thom(G-v,H)$, whereas
  $\thom(H,G)$ collapses onto $\thom(H,G-v)$. 
The maps $i_H$ and $f^H$ are strong deformation retractions.
\end{thm}

\begin{figure}[hbt]
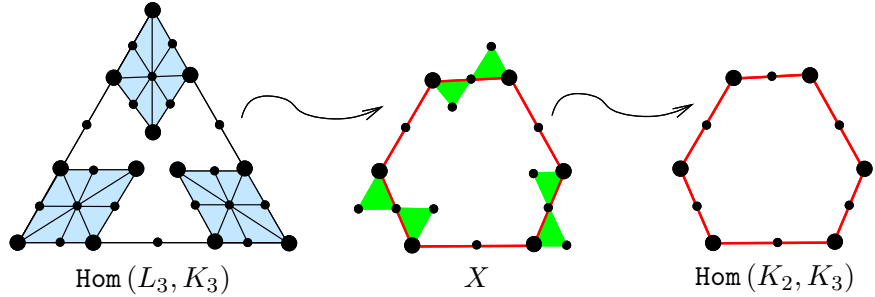

\begin{center}
  \begin{picture}(0,0)%
    \includegraphics{coll1.pstex}%
  \end{picture}%
  \input{coll1.pstex_t}%
  
\end{center}
\caption{A two-step folding of the first argument in $\thom(L_3,K_3)$.}
\label{fig:coll1}
\end{figure}

\pr First we show that $\bd\thom(G,H)$ collapses onto
$\bd\thom(G-v,H)$. Identify $\cp(\thom(G-v,H))$ with the subposet of
$\cp(\thom(G,H))$ consisting of all $\eta$, such that
$\eta(v)=\eta(u)$. Let $X$ be the subposet consisting of all
$\eta\in\cp(\thom(G,H))$ satisfying $\eta(v)\supseteq\eta(u)$. Then
$\cp(\thom(G-v,H))\subseteq X \subseteq\cp(\thom(G,H))$. Consider
order preserving maps
\[\cp(\thom(G,H))\stackrel{\alpha}{\longrightarrow} X\stackrel{\beta}{\longrightarrow}\cp(\thom(G-v,H)),\]
defined by
\[\alpha\eta(x)=\begin{cases}
        \eta(u)\cup\eta(v), & \text{for }x=v;\\
        \eta(x), & \text{otherwise;}
\end{cases}
\qquad\quad
\beta\eta(x)=\begin{cases}
        \eta(u), & \text{for }x=v;\\
        \eta(x), & \text{otherwise;}
\end{cases}
\]
for all $x\in V(G)$; see Figure \ref{fig:coll1}. Maps $\alpha$ and
$\beta$ are well-defined because $G-v$ is a~fold of~$G$. Clearly
$\beta\circ\alpha=\cp(i_H)$, $\alpha$ is an~ascending closure
operator, and $\beta$ is a~descending closure operator. Since
Im$\,\cp(i_H)=\cp(\thom(G-v,H))$, the statement follows from
Theorem~\ref{thm1}.

We show that $\thom(H,G)$ collapses onto $\thom(H,G-v)$ by presenting
a~sequence of elementary collapses. Denote $V(H)=\{x_1,\dots,x_t\}$.
For $\eta\in\cp(\thom(H,G))$, let $1\leq i(\eta)\leq t$ be the minimal
index such that $v\in\eta(x_{i(\eta)})$. Write $\cp(\thom(H,G))$ as
a~disjoint union $A\cup B\cup\cp(\thom(H,G-v))$, defined as follows:
for $\eta\in A\cup B$ we have $\eta\in A$ if
$u\notin\eta(x_{i(\eta)})$, and we have $\eta\in B$ otherwise.

There is a~bijection $\varphi:A\ra B$ which adds $u$ to
$\eta(x_{i(\eta)})$ without changing other values of $\eta$. Adding
$u$ to $\eta(x_{i(\eta)})$ yields an~element in $\cp(\thom(H,G))$
since $G-v$ is a~fold of~$G$. Clearly, $\varphi(\alpha)$ covers
$\alpha$, for all $\alpha\in A$. We take the set
$\{(\alpha,\varphi(\alpha))\,|\,\alpha\in A\}$ to be our collection of
the elementary collapses. These are ordered lexicographically after
the pairs of integers $(i(\alpha),-\dim\alpha)$.

Let us see that these collapses can be performed in this lexicographic
order. Take $\eta>\alpha$, $\eta\neq\varphi(\alpha)$.  Assume
$i(\eta)=i(\alpha)$.  If $\eta\in B$, then
$\eta=\varphi(\tilde\alpha)$, $i(\tilde\alpha)=i(\alpha)$, and
$\dim\tilde\alpha>\dim\alpha$.  Otherwise $\eta\in A$ and
$\dim\eta>\dim\alpha$. The third possibility is that
$i(\eta)<i(\alpha)$. In either case $\eta$ has been removed
before~$\alpha$.  \qed

\vskip5pt

Instead of verifying that the sequence of collapses is correct in the
last paragraph of the proof we could simply notice that the defined
matching is acyclic and derive the result by discrete Morse
theory,~\cite{Fo}.
   
\begin{rem} 
In analogy with the first part of the proof, we can show
that $\bd\thom(H,G)$ collapses onto $\bd\thom(H,G-v)$ by rewriting
$\cp(f^H)$ as a~composition of two closure operators.
\end{rem} 

Indeed, let $Y$ be the
subposet consisting of all $\eta\in\cp(\thom(H,G))$ such that for all
$x\in V(H)$, $\eta(x)\cap\{u,v\}\neq\{v\}$, i.e., for any $x\in V(H)$
we have: if $v\in\eta(x)$, then $u\in\eta(x)$. Then $\cp(\thom(H,G-v))
\subseteq Y \subseteq\cp(\thom(H,G))$. Consider order preserving maps
\[\cp(\thom(H,G))\stackrel{\varphi}{\longrightarrow} Y
\stackrel{\psi}{\longrightarrow}\cp(\thom(H,G-v)),\]
defined by
\[\varphi\eta(x)=\begin{cases}
        \eta(x)\cup\{u\}, & \text{if }v\in\eta(x);\\
        \eta(x), & \text{otherwise;}
\end{cases}
\qquad
\psi\eta(x)=\begin{cases}
        \eta(x)\setminus\{v\}, & \text{if }v\in\eta(x);\\
        \eta(x), & \text{otherwise;}
\end{cases}
\]
for all $x\in V(H)$. 

The map $\phi$ is well-defined because $G-v$ is a~fold of~$G$, and the
map $\psi$ is well-defined by the construction of~$Y$.  We see that
$\psi\circ\varphi=\cp(f^H)$, $\varphi$ is an~ascending closure
operator, and $\psi$ is a~descending closure operator. Since
Im$\,\cp(f_H)=\cp(\thom(H,G-v))$, the statement follows from
Theorem~\ref{thm1}.

\vskip10pt

\nin {\bf Acknowledgments.} We would like to thank P\'eter Csorba for 
useful comments and the Swiss National Science Foundation and
ETH-Z\"urich for the financial support of this research.

\end{document}